\documentclass[12pt]{article}

\usepackage{amssymb,hyperref}

\usepackage{latexsym}
\newtheorem{theoreme}{{\bf Th\'eor\`eme}}[section]

\newtheorem{corollaire principal}[principal]{\bf Corollaire}
\newtheorem{proposition}[theoreme]{{\bf Proposition}}
\newtheorem{lemme}[theoreme]{{\bf Lemme}}
\newtheorem{corollaire}[theoreme]{{\bf Corollaire}}
\newtheorem{definition}[theoreme]{{\bf D\'efinition}}

\newenvironment{demonstration}{\noindent{\bf D\'emonstration
}}{\nolinebreak $\Box $\hspace{-2.15mm}\raisebox{1.25mm}{.} \medskip}

\newenvironment{demonstration du lemme}{\noindent{\bf D\'emonstration du lemme
}}{\nolinebreak $\Box $\hspace{-2.15mm}\raisebox{1.25mm}{.} \medskip}

\def\RR{{\bf R}}
\def\CC{{\bf C}}
\def\ZZ{{\bf Z}}

\begin{document}

\title{{\Large \bf  Connexions affines et projectives sur les surfaces complexes compactes}}

\author{{\normalsize
{\bf Sorin DUMITRESCU}}}

\date{Mai 2008}

\maketitle

\vspace{0.3cm}

\noindent{\bf{ R\'esum\'e.}} 
Soit $(S ,\nabla)$ une surface complexe compacte connexe  munie d'une connexion affine holomorphe sans torsion. Nous d\'emontrons que  $\nabla$ est localement model\'ee sur une connexion affine  invariante par translations sur $\CC^2$ (en particulier, $\nabla$ est localement homog\`ene),  sauf si  $S$ est un fibr\'e elliptique principal au-dessus d'une surface de genre $g \geq 2$, de premier nombre de Betti impair et $\nabla$ est une connexion affine holomorphe sans torsion g\'en\'erique sur  $S$, auquel cas l'alg\`ebre de Lie des champs de Killing locaux est de dimension un, engendr\'ee par le champ fondamental de la fibration principale.\\
Nous en d\'eduisons que toute connexion projective holomorphe normale sur une surface complexe compacte est  plate.

\section{Introduction}

Suite \`a des  nombreux travaux portant  sur  les surfaces  complexes  qui poss\`edent des connexions affines holomorphes~\cite{Gunning,
Vitter, Suwa, Maehara}, Inoue, Kobayashi et Ochiai   \'etablissent dans~\cite{IKO}  la liste des surfaces complexes (connexes) compactes qui admettent des connexions affines holomorphes
(non n\'ecessairement plates). Une telle surface est biholomorphe (\`a rev\^etement fini pr\`es) \`a un tore complexe, une surface de Kodaira primaire, une surface de Hopf affine,
une surface d'Inoue ou \`a un fibr\'e principal en courbes elliptiques sur une surface de Riemann de genre $g \geq 2$, de premier nombre de Betti impair. Par ailleurs, il est
montr\'e dans~\cite{IKO} que toutes ces surfaces poss\`edent des connexions affines holomorphes sans torsion {\it plates}. Ces connexions sont localement 
model\'ees sur $\CC^2$ muni de sa connexion standard et, par cons\'equent, toutes ces surfaces poss\`edent au moins une {\it structure affine holomorphe}, i.e. un atlas 
\`a valeurs dans des ouverts de $\CC^2$ et \`a changements de carte localement constants dans le groupe affine $GL(2, \CC) \ltimes \CC^2$.

La classification des surfaces complexes compactes qui poss\`edent des connexions projectives holomorphes a \'et\'e obtenue dans~\cite{KO1,KO2}. Ces surfaces sont 
dans la liste suivante : $P^2(\CC)$, les quotients non ramifi\'es de l'espace hyperbolique complexe $H^2_{\CC}$ ou les surfaces admettant une connexion affine holomorphe.
Toutes ces surfaces poss\`edent des connexions projectives holomorphes {\it plates}, i.e. localement model\'ees sur $P^2(\CC)$, et donc des {\it structures projectives holomorphes} (voir section 2).

Dans~\cite{Klingler} Klingler classifie  les structures affines holomorphes et les structures projectives holomorphes sur les surfaces complexes compactes.

Le but de cet article est d'\'etudier la g\'eom\'etrie locale  des connexions affines et projectives holomorphes  sur les surfaces complexes compactes. Nos r\'esultats sont les suivants.

\begin{theoreme} \label{principal} Soit $(S ,\nabla)$ une surface complexe compacte connexe munie d'une connexion affine holomorphe sans torsion.

i) Si $S$ n'est pas biholomorphe \`a un fibr\'e principal elliptique au-dessus d'une surface de Riemann de genre $g \geq 2$, de premier nombre de Betti impair, alors $\nabla$ est localement model\'ee sur
une connexion affine invariante par translations sur $\CC^2$. En particulier, $\nabla$ est  localement homog\`ene.

ii) Supposons que  $S$ est un fibr\'e principal holomorphe en courbes elliptiques au-dessus d'une surface $\Sigma$ de genre $g \geq 2$, de premier nombre de Betti impair. L'espace des connexions affines holomorphes sans torsion  sur $S$ est biholomorphe \`a  $$ (H^0(\Sigma, K_{\Sigma}) )^2 \times H^0(\Sigma, K_{\Sigma}^{\otimes 2}) \simeq \CC^{5g-3},$$
o\`u $K_{\Sigma}$ d\'esigne le fibr\'e canonique de $\Sigma$.

Toutes ces connexions sont invariantes par la fibration principale (i.e. admettent le champ fondamental de la fibration principale comme champ de Killing). L'espace des connexions
plates est de codimension complexe $g$. Les  connexions affines holomorphes sans torsion non plates  sont non  localement homog\`enes : l'alg\`ebre des champs de  Killing locaux est de dimension un, engendr\'ee par le champ fondamental de la  fibration principale.  
\end{theoreme}

 Les connexions construites au  point (ii) du th\'eor\`eme pr\'ec\'edent fournissent les premiers exemples de $G$-structures holomorphes rigides  (voir section~\ref{definitions}) sur une vari\'et\'e compacte
qui ne sont pas localement homog\`enes (voir~\cite{D1,D2,HM,McKay,McKay2} pour des r\'esultats d'homog\'en\'eit\'e locale et de rigidit\'e concernant ces structures).

La classification du th\'eor\`eme~\ref{principal} nous permet d'obtenir le 

\begin{corollaire} \label{cconnexions} Soit $(S, \nabla)$ une surface complexe compacte munie d'une connexion affine holomorphe sans torsion. Alors $\nabla$ est projectivement plate.
\end{corollaire}

En utilisant le corollaire~\ref{cconnexions}, nous obtenons  le r\'esultat  suivant  :

\begin{theoreme}  \label{projectif} Toute connexion projective holomorphe normale sur une surface complexe compacte est  plate.
\end{theoreme}

Il convient de mentionner ici que  les tores de  dimension strictement sup\'erieure \`a deux poss\`edent  des connexions affines holomorphes sans torsion  (invariantes par translations) 
qui ne sont pas projectivement plates (voir, par exemple, ~\cite{McKay3}).

\`A titre de motivation rappelons  la conjecture de   B. McKay~\cite{McKay} qui affirme que sur les vari\'et\'es complexes compactes, les connexions de Cartan paraboliques holomorphes  sont  {\it localement homog\`enes}. Le th\'eor\`eme~\ref{projectif} repond affirmativement \`a la question en dimension complexe deux. En effet, en dimension deux,  les connexions  de Cartan paraboliques holomorphes sont
essentiellement des connexions projectives holomorphes normales ou  des structures conformes holomorphes~\cite{McKay}. Ces derni\`eres sont toujours  plates,  d'apr\`es le  r\'esultat classique, d\^u \`a Gau{\ss}, d'existence de coordonn\'ees isothermes.
Pour la  classification des surfaces complexes compactes admettant   des structures conformes holomorphes  le lecteur pourra consulter~\cite{KO}.

Ajoutons que si  la classification des vari\'et\'es complexes compactes de dimension $3$ poss\`edant des connexions projectives holomorphes normales est connue dans le cadre alg\'ebrique~\cite{JR}, elle  reste ouverte    en g\'en\'eral.

Le plan de l'article est le suivant. Dans la section~\ref{definitions} nous pr\'ecisons  les d\'efinitions et d\'emontrons le  th\'eor\`eme~\ref{principal} dans le cas o\`u la surface $S$ n'est pas une surface d'Inoue ou un fibr\'e principal elliptique  au-dessus d'une surface de Riemann de genre $g \geq 2$, de premier nombre de Betti impair. 

Dans la section 3 nous traitons le cas des surfaces d'Inoue. Au passage nous d\'emontrons que sur ces  surfaces, toutes les structure g\'eom\'etriques holomorphes  sont localement homog\`enes.

Dans la section 4 nous r\`eglons le cas restant des fibr\'es principaux elliptiques et  nous ach\'evons la preuve  du
th\'eor\`eme~\ref{principal} et de son corollaire~\ref{cconnexions}. 

Enfin \`a la  section~\ref{Connexions projectives}  nous d\'emontrons le th\'eor\`eme~\ref{projectif}. 

{\it Je tiens \`a remercier  B. McKay  pour des conversations utiles et agr\'eables sur le sujet. Je remercie  \'egalement  G. Dloussky.

 Ce travail a \'et\'e partiellement soutenu par l'ANR Symplexe BLAN 06-3-137237.}

\section{D\'efinitions. Contexte}  \label{definitions}

Le livre de Gunning~\cite{Gunning} fournit  la premi\`ere \'etude syst\'ematique des connexions affines et projectives holomorphes. Contrairement au cas r\'eel, il est montr\'e dans~\cite{Gunning} qu'une vari\'et\'e complexe compacte de dimension $n \geq 2$  n'admet pas toujours une connexion affine ou projective holomorphe et que la
pr\'esence d'une telle connexion entra\^ine l'annulation de certaines classes caract\'eristiques.

Moralement, une connexion affine holomorphe $\nabla$ permet de d\'eriver les champs de vecteurs holomorphes  (locaux)  $Y$, le long des orbites d'un deuxi\`eme champ
de vecteurs holomorphe (local) $X$ et d'obtenir comme r\'esultat un champ de vecteurs holomorphe (local) $\nabla_{X}Y$. 

Si l'on multiplie $X$ et $Y$
par des fonctions (locales) $f$, cette op\'eration est tensorielle en $X$ et satisfait \`a  la r\`egle  $\nabla_{X}fY=(X \cdot f)Y+f\nabla_{X}Y$. Il d\'ecoule  que la diff\'erence entre deux connexions affines holomorphes  sur une vari\'et\'e $S$
est un tenseur, section holomorphe du fibr\'e vectoriel $T^*S \otimes T^*S \otimes TS$. 

Si, de plus, $\nabla_{X}Y=\nabla_{Y}X$, pour tous les champs de vecteurs  $X$ et $Y$ qui commutent, alors $\nabla$ est dite {\it sans torsion}.

Dans la suite nous donnerons une d\'efinition \'equivalente qui permet de voir les connexions affines et les connexions projectives comme des cas particuliers de $G$-structures~\cite{Kobayashi}.

Consid\'erons une vari\'et\'e $S$ de dimension $n$ (pour  ce qui est de cet article,  on pourra se contenter de penser que $n=2$).

Rappelons que pour tout entier positif $r$, le groupe $D^r$, des $r$-jets en $0$ de germes de biholomorphismes locaux  de $\CC^n$ qui fixent $0$ est un groupe lin\'eaire alg\'ebrique qui co\"{\i}ncide avec $GL(n, \CC)$, pour $r=1$, et avec une extension de $GL(n,\CC)$ par le groupe additif des formes bilin\'eaires sym\'etriques sur $\CC^n$, si $r=2$.

Le fibr\'e des $r$-rep\`eres $R^r(M)$ de $M$, autrement dit le fibr\'e des $r$-jets   en $0$ de germes de biholomorphismes locaux entre $\CC^n$ et $M$ est un fibr\'e principal au-dessus de $M$ de groupe structural $D^r$. Nous suivons~\cite{DG,Gro} et donnons la 

\begin{definition} Une structure g\'eom\'etrique holomorphe  $\tau$ (d'ordre $r$) sur $M$ est une application holomorphe, $D^r$-\'equivariante, $\tau : R^r(M) \to Z$,  avec  $Z$  une vari\'et\'e alg\'ebrique 
munie d'une action alg\'ebrique de $D^r$.
\end{definition}

L'application $\tau$ s'interpr\`ete comme une section holomorphe du fibr\'e de fibre $Z$, associ\'e au fibr\'e principal $R^r(M)$ via l'action de $D^r$ sur $Z$.

Un biholomorphisme local $f$ de $M$ agit naturellement sur les sections du fibr\'e pr\'ec\'edent. Si cette action pr\'eserve $\tau$, alors $f$ est une {\it isom\'etrie locale}  de $\tau$.
Si les isom\'etries locales agissent transitivement sur $M$, alors la structure g\'eom\'etrique $\tau$ est dite {\it localement homog\`ene.}

Si l'image de $\tau$ dans $Z$ est exactement une $D^r$-orbite, qui s'identifie alors  \`a un espace homog\`ene $D^r/G$, o\`u $G$ est le sous-groupe de $D^r$ qui stabilise un point de l'image, alors $\tau$ s'interpr\`ete comme une section d'un fibr\'e de fibre $D^r/G$. Cette section fournit une r\'eduction du groupe structural de $R^r(M)$ au sous-groupe $G$. Une telle
structure g\'eom\'etrique est dite une $G$-structure holomorphe~\cite{Kobayashi}.

Consid\'erons le sous-groupe $G$ de $D^2$,  isomorphe \`a $D^1$,  constitu\'e par les $2$-jets en $0$ d'isomorphismes  lin\'eaires de $\CC^n$, ainsi que le sous-groupe $H$ constitu\'e par les $2$-jets en $0$ de transformations projectives de $P^n(\CC)$ qui fixent $0$. La d\'efinition \'equivalente suivante des connexions affines et projectives, vues  comme $G$-structures est classique~\cite{Kobayashi Nagano, MM}:

\begin{definition} Une connexion affine (holomorphe) sans torsion sur $M$ est une r\'eduction (holomorphe) du groupe structural de $R^2(M)$ \`a $G$.

                                  Une connexion projective (holomorphe) normale sur $M$ est une r\'eduction (holomorphe) du groupe structural de $R^2(M)$ \`a $H$.
\end{definition}

L'inclusion canonique de $G$ dans $H$ permet d'associer \`a chaque connexion affine sans torsion une unique connexion projective normale.

Le lecteur pourra \'egalement se r\'ef\'erer \`a~\cite{MM}, pour voir que la d\'efinition pr\'ec\'edente co\"{\i}ncide avec celle adopt\'ee dans~\cite{KO1,KO2}. Une connexion projective normale d\'efinit une unique connexion  de Cartan projective~\cite{MM}.

Une connexion affine sans torsion (projective normale) est dite {\it plate} si elle est localement isomorphe \`a la connexion affine (projective) standard de l'espace $\CC^n$ ( de l'espace projectif $P^n(\CC)$). Une connexion affine sans torsion est plate si et seulement si son tenseur de courbure s'annule~\cite{Gunning}. La platitude des connexions projectives est
test\'ee par l'annulation du tenseur courbure de Weyl (voir, par exemple,~\cite{McKay3}).

Comme dans le cas affine, une connexion projective normale plate sur $M$ donne naissance \`a une {\it structure  projective}. Autrement dit, il existe un atlas de $M$ \`a valeurs dans
des ouverts de $P^n(\CC)$ tel que les applications de changement de carte soient des transformations projectives. Dans ce cas, il existe une {\it application d\'eveloppante} d\'efinie
sur le rev\^etement universel $\tilde M$ \`a valeurs dans $P^n(\CC)$ et un {\it morphisme d'holonomie} du groupe fondamental de $M$ dans le groupe $PSL(n+1, \CC)$  qui rend l'application d\'eveloppante \'equivariante~\cite{Gunning}.

Rappelons que les connexions affines et projectives sont des structures {\it rigides} au sens de Gromov~\cite{DG,Gro}, i.e. le pseudo-groupe des isom\'etries locales est un pseudo-groupe de Lie, engendr\'e par une alg\`ebre de Lie de dimension finie appel\'ee {\it l'alg\`ebre des champs de Killing locaux}. Dans le cas analytique la dimension de cette alg\`ebre ne d\'epend 
pas du point (de $M$)~\cite{Amores,DG,Gro}. 

Une isom\'etrie (locale) entre deux connexions projectives normales sera \'egalement appel\'ee {\it \'equivalence projective} (locale).

Revenons \`a pr\'esent aux  connexions affines holomorphes sur les surfaces complexes.

{\it Connexions affines invariantes par translations sur $\CC^2$.}

Les cas les plus simples de connexions affines holomorphes sont celles qui sont invariantes par transations sur $\CC^2$. Autrement dit, dans les coordonn\'ees
 canoniques $(z_{1},z_{2})$, on d\'efinit $\nabla_{\frac{\partial}{\partial z_{i}}}\frac{\partial}{\partial z_{j}}=\Gamma_{ij}^1 \frac{\partial}{\partial z_{1}} +\Gamma_{ij}^2 \frac{\partial}{\partial z_{2}}$, avec
 $\Gamma_{ij}^1$ et $\Gamma_{ij}^2$ des fonctions constantes, quelque soient $1 \leq i,j \leq 2$.
 
 Ces connexions descendent  sur tout tore complexe de dimension $2$, quotient de $\CC^2$ par un r\'eseau de translations. Par ailleurs, comme le fibr\'e tangent au tore est holomorphiquement
 trivial, toute connexion affine holomorphe sur un tore  s'obtient par ce proc\'ed\'e et est invariante par translations, donc (localement) homog\`ene.

 La connexion est sans torsion si et seulement si
 $\Gamma_{ij}^1=\Gamma_{ji}^1$ et $\Gamma_{ij}^2=\Gamma_{ji}^2$.
 
 Un calcul direct montre que la courbure est nulle  si et seulement si les coefficients  $\Gamma_{ij}^1$ et $\Gamma_{ij}^2$ satisfont quatre relations quadratiques (voir, par exemple, ~\cite{IKO}). En particulier, les tores de dimension \'egale (ou sup\'erieure) \`a deux poss\`edent des connexions affines holomorphes sans torsion non plates.

En rassemblant les r\'esultats existant dans la litt\'erature~\cite{IKO,Vitter,Wall} on obtient la proposition suivante :

\begin{proposition} \label{debut}  Soit $\nabla$ une connexion affine holomorphe sur  une surface complexe compacte $S$ non biholomorphe \`a un fibr\'e principal elliptique au-dessus d'une surface de Riemann de genre $g \geq 2$, de premier nombre de Betti impair,  ou \`a une surface d'Inoue, alors  $\nabla$ est localement isom\'etrique \`a une connexion affine  invariante par translations sur $\CC^2$.
\end{proposition}

En effet, le r\'esultat de~\cite{IKO} implique  les   possibilit\'ees  suivantes :

-$S$ est un tore complexe. Dans ce cas nous venons de pr\'esenter la preuve.

-$S$ est une surface de Hopf. Dans ce cas la courbure est n\'ecessairement nulle~\cite{IKO} et si $\nabla$ est sans torsion, elle est localement isomorphe \`a la connexion standard
sur $\CC^2$.

-$S$ est une surface de Kodaira primaire (i.e. un fibr\'e en courbes elliptiques au-dessus d'une courbe elliptique).  Les connexions
affines sur ce type de surfaces ont \'et\'e \'etudi\'ees par Vitter. Dans~\cite{Vitter} (page 239) Vitter  montre que ces connexions sont n\'ecessairement localement isom\'etriques \`a des connexions invariantes par translations sur $\CC^2$.

\section{Surfaces d'Inoue}

Les  surfaces d'Inoue sont des surfaces complexes compactes de type $VII_{0}$, de dimension alg\'ebrique nulle (i.e. sans fonctions m\'eromorphes non constantes) et avec  $b_{2}=0$~\cite{Barth-Peters}. Ces surfaces ont \'et\'e \'etudi\'ees dans~\cite{Inoue}, o\`u il est montr\'e que chaque telle surface est obtenue comme quotient de $H \times \CC$ par un groupe $\Gamma$ de
transformations affines de $\CC^2$ pr\'eservant cet ouvert ($H$ \'etant le demi-plan sup\'erieur de Poincar\'e). En particulier, il d\'ecoule de~\cite{Inoue} que chaque surface d'Inoue h\'erite d'une structure affine induite
par la structure affine canonique de $\CC^2$.  Il est montr\'e dans~\cite{Klingler} (lemme 4.3) que cette structure affine est unique.
Nous   d\'emontrons ici le r\'esultat plus g\'en\'eral suivant  :

\begin{proposition}  \label{Inoue}   Une surface d'Inoue admet une unique connexion affine holomorphe sans torsion. 
\end{proposition}

 L'unique connexion affine de l'\'enonc\'e pr\'ec\'edent est n\'ecessairement  plate. Les propositions~\ref{Inoue} et~\ref{debut} ach\`event la preuve du point (i) du th\'eor\`eme~\ref{principal}.
 
 Avant de passer \`a la preuve de la proposition~\ref{Inoue} rappelons qu'\`a rev\^etement double non ramifi\'e pr\`es,  les surfaces d'Inoue s'obtiennent par l'un des deux proc\'ed\'es suivants~\cite{Inoue}.
 
 {1. \it Surfaces $S_{M}$.} Consid\'erons une matrice $M \in SL(3,\ZZ)$ avec valeurs propres $\alpha, \beta, \bar{\beta}$ tels que $\alpha >1$ et $\beta \neq \bar{\beta}$. Choisissons
 $(a_{1},a_{2},a_{3})$ un vecteur propre r\'eel associ\'e \`a la valeur propre $\alpha$ et $(b_{1},b_{2},b_{3})$ un vecteur propre associ\'e \`a la valeur propre $\beta$.
 
 Consid\'erons le groupe  $\Gamma$ de transformations (affines)   de $\CC^2$ engendr\'e par  :

$\gamma_{0}(w,z)=(\alpha w, \beta z),$

$\gamma_{i}(w,z)=(w+a_{i}, z+b_{i}),$ avec $i=1,2,3.$

L'action de $\Gamma$ sur $\CC^2$ pr\'eserve $H \times \CC$ et le quotient est la surface complexe compacte $S_{M}$.

{2. \it Surfaces $S^+_{N,p,q,r,t}$.}  Soit $N =(n_{ij})  \in SL(2, \ZZ)$ une matrice diagonalisable sur $\RR$ avec valeurs propres $\alpha >1$ et $ \alpha^{-1}$  et vecteurs propres $(a_{1},a_{2})$ et
$(b_{1},b_{2})$ respectivement. Choisissons $r \in \ZZ^*, p,q  \in \ZZ , t \in \CC$ et des nombres r\'eelles $c_{1}, c_{2}$ solutions de l'\'equation 

$$(c_{1},c_{2})=(c_{1},c_{2}) N^t +(e_{1},e_{2}) + \frac{1}{r}(b_{1}a_{2}-b_{2}a_{1})(p,q),$$

o\`u $e_{i}= \frac{1}{2}n_{i1}(n_{i1}-1)a_{1}b_{1}  + \frac{1}{2}n_{i2}(n_{i2}-1)a_{2}b_{2}  +n_{i1}n_{i2}b_{1}a_{2}$ et $N^t$ d\'esigne la transpos\'ee de $N$.

Dans ce cas $\Gamma$ est engendr\'e par les transformations 

$\gamma_{0}(w,z)=(\alpha w,z+t),$

$\gamma_{i}(w,z)=(w+a_{i}, z+b_{i}w +c_{i}),~~(i=1,2),$

$\gamma_{3}(w,z)=(w, z+ r^{-1}(b_{1}a_{2} -b_{2}a_{1})).$

Ce groupe est discret et agit proprement et discontinument sur $H \times \CC$ admettant comme quotient la surface complexe compacte $S^+_{N,p,q,r,t}$.

Dans les deux cas  notons $\Lambda = \lbrack \Gamma,  \Gamma \rbrack$ le commutateur de $\Gamma$. Le groupe $\Lambda$ est engendr\'e par $\gamma_1, \gamma_{2}$ et $\gamma_{3}$.
 
 Le lemme suivant (\cite{Inoue}, lemma 3, lemma 4) sera essentiel :
 
 \begin{lemme} \label{constant}  Les fonctions holomorphes $\Lambda$-invariantes sur $H \times \CC$ sont n\'ecessairement constantes.
 \end{lemme}
 
 Passons \`a pr\'esent \`a la preuve de la proposition~\ref{Inoue}.
 
 \begin{demonstration}  Consid\'erons la surface d'Inoue $S$ munie de sa connexion affine holomorphe  sans torsion plate  $\nabla_{0}$, h\'erit\'ee par la connexion canonique de $H \times \CC$.
 La diff\'erence entre une  connexion holomorphe affine arbitraire $\nabla$  sur $S$ et $\nabla_{0}$ est une section holomorphe $\omega$ du fibr\'e vectoriel 
 $T^*S \otimes T^*S \otimes TS$. L'image r\'eciproque de ce tenseur sur le rev\^etement  $H \times \CC$ de $S$ est une section du m\^eme type $\tilde{\omega}$ invariante
par le groupe de rev\^etement $\Gamma$.

En les coordonn\'ees globales   $(w, z) \in H  \times \CC$ du rev\^etement universel,  l'expression de $\tilde{\omega}$ sera  :
$$ f_{11}(w, z) dz \otimes dz \otimes \frac{\partial}{\partial z} + 
   f_{12}(w, z) dz \otimes dw  \otimes \frac{\partial}{\partial z}+
    f_{21}(w, z) dz \otimes dz \otimes \frac{\partial}{\partial w}+
     f_{22}(w, z) dz \otimes dw  \otimes \frac{\partial}{\partial w}+$$
    $$  g_{11}(w, z) dw  \otimes dz \otimes \frac{\partial}{\partial z}+
       g_{12}(w, z) dw  \otimes dw \otimes \frac{\partial}{\partial z}+
        g_{21}(w, z) dw  \otimes dz \otimes \frac{\partial}{\partial  w}+
         g_{22}(w, z) dw  \otimes dw  \otimes \frac{\partial}{\partial w},$$
         
         avec $f_{ij},g_{ij}$ des fonctions holomorphes quelque soient $1 \leq i,j \leq 2$.

Nous consid\'erons seulement des  connexions $\nabla$ sans torsion, ce qui conduit \`a : $f_{12}=g_{11}$ et  $f_{22}=g_{21}$.

{\it Commen\c cons par le cas des surfaces $S_{M}.$ }

L'invariance de $\tilde {\omega}$ par le groupe $\Lambda$ engendr\'e par les translations $\gamma_{1}, \gamma_{2}$ et $\gamma_{3}$, implique que
les fonctions $f_{ij}$ et $g_{ij}$ sont $\Lambda$-invariantes et, gr\^ace au lemme~\ref{constant}, constantes.

L'invariance de $\tilde{\omega}$ par $\gamma_{0}$ conduit \`a $f_{11}= \beta f_{11}, f_{12}=\alpha f_{12}, f_{21}=\frac{\beta^2}{\alpha} f_{21}, f_{22}= \beta f_{22}, g_{21}=\beta g_{21}$
et  $g_{22}= \alpha g_{22}$. Comme les nombres complexes $\alpha, \beta$ et $\frac{\beta^2}{\alpha}$ sont $\neq 1$,  toutes les fonctions $f_{ij}, g_{ij}$ sont nulles et, par cons\'equent,
$\tilde{\omega}=0$.

Il vient que la connexion $\nabla_{0}$ est l'unique connexion affine sans torsion sur $S_{M}$.

{\it On passe aux cas des surfaces $S^+_{N,p,q,r,t}$. }

L'invariance de $\tilde{\omega}$ par la translation $\gamma_{3}$ implique l'invariance des fonctions $f_{ij},g_{ij}$ par $\gamma_{3}$.

L'invariance de $\tilde{\omega}$ sous l'action de $\gamma_{1}$ et $\gamma_{2}$ s'exprime sous la forme suivante :

 $(1)   f_{11}(w,z)=f_{11}(\gamma_{i} (w,z)) -b_{i} f_{21}(\gamma_{i}(w,z))$
                  
                $  (2)     f_{12}(w,z)=f_{12}(\gamma_{i}(w,z)) -2b_{i}^2f_{21}(\gamma_{i}(w,z))  -b_{i} f_{22} (\gamma_{i}(w,z)) +2b_{i} f_{11}(\gamma_{i}(w,z))$
                  
                  $(3)     f_{21}(w,z)=f_{21}(\gamma_{i}(w,z))$
                  
                 $ (4)     f_{22}(w,z)=2b_{i}f_{21}(\gamma_{i}(w,z))  +f_{22}(\gamma_{i} (w,z))$
                  
                 $ (5)     g_{12}(w,z)= g_{12}(\gamma_{i} (w,z)) +b_{i}^2 f_{11}(\gamma_{i}(w,z))   +  b_{i} f_{12}(\gamma_{i}(w,z))-b_{i}^3f_{21}(\gamma_{i} (w,z)) -b_{i}^2f_{22}(\gamma_{i}(w,z)) -b_{i}g_{22}(\gamma_{i}(w,z))$
                  
                  $(6) g_{22}(w,z)=g_{22}(\gamma_{i}(w,z)) +b_{i}f_{22}(\gamma_{i}(w,z))+b_{i}^2f_{21}(\gamma_{i}(w,z)),$

                  pour $i=1,2$.
                  
          L'\'equation $(3)$ montre que la fonction $f_{21}$ est invariante par $\gamma_{i}$, pour $i=1,2$. Comme on a vu que cette fonction est \'egalement $\gamma_{3}$-invariante,
          elle est $\Lambda$-invariante et, par le lemme~\ref{constant}, $f_{21}$ est constante.
          
          Or, l'invariance de $\tilde{\omega}$ par $\gamma_{0}$ implique alors $f_{21}= \frac{1}{\alpha}f_{21}$ et, comme $\alpha \neq 1$, $f_{21}=0$.
          
          La premi\`ere \'equation montre alors que la fonction $f_{11}$ est invariante par $\gamma_{1}$ et $\gamma_{2}$. Il vient que  $f_{11}$ est $\Lambda$-invariante et donc  constante.
          
          Aussi la quatri\`eme \'equation implique que $f_{22}$ est invariante par $\gamma_{1}$ et $\gamma_{2}$ et donc $\Lambda$-invariante et, par le lemme~\ref{constant}, constante. 
          
          En d\'erivant par rapport \`a la variable $z$ l'\'equation $(6)$ on obtient alors que la fonction 
          $\displaystyle \frac{\partial g_{22}}{\partial z}$ est $\gamma_{i}$-invariante pour $i=1,2$. Comme cette fonction holomorphe est \'egalement $\gamma_{3}$-invariante, le lemme~\ref{constant}
          implique que $\displaystyle \frac{\partial g_{22}}{\partial z}$ est une constante $C \in \CC$ et, par cons\'equent, $g_{22}(w,z)=Cz+g(w)$, avec $g$ fonction holomorphe de la variable $w$.
          
          Par ailleurs, l'invariance de $\tilde{\omega}$ par $\gamma_{0}$ implique que $g_{22}(w,z)=\alpha g_{22}(\alpha w,z)$. Comme $\alpha \neq 1$, ceci implique $C=0$
          et donc  $g_{22}(w,z)=g(w)=\alpha g(\alpha w)$.
          
          En revenant \`a l'\'equation $(6)$ on a $g(w)-g(w+a_{i})=b_{i}f_{22}$. L'homog\'en\'eit\'e de $g$ implique alors que $\displaystyle g(w)-g(w + \frac{a_{i}}{\alpha^k})=\alpha^kb_{i}f_{22}$, pour
          tout entier positif $k$. On fait tendre $k$ vers l'infini et,  comme $\alpha > 1$, la continuit\'e de $g$ en $w$ implique que $f_{22}=0$.
          
          L'\'equation $(6)$ implique alors que $g_{22}$ est $\gamma_{i}$-invariante, donc constante, ainsi que $g$.
         De plus,  l'homog\'en\'eit\'e de $g$ implique que $g$ est nulle. Il vient que $g_{22}=0$.
          
          La deuxi\`eme \'equation devient alors $ f_{12}(w,z)=f_{12}(\gamma_{i}(w,z)) +2b_{i} f_{11}(\gamma_{i}(w,z))$. On obtient comme avant que
          $\displaystyle \frac{\partial f_{12}}{\partial z}$ est une constante et ensuite que $f_{12}=f_{11}=0$.
          
          L'\'equation $(5)$ implique alors que $g_{12}$ est $\gamma_{i}$-invariante pour $i=1,2$. Il vient que $g_{12}$ est constante. L'invariance par $\gamma_{0}$ de
         $\tilde{\omega}$ implique alors $g_{12}= \alpha^2g_{12}$ et donc $g_{12}=0$.
         
         On vient de d\'emontrer que $\tilde{\omega}=0$, ce qui ach\`eve la preuve du second  cas.
 \end{demonstration}

Nous d\'emontrons ici le th\'eor\`eme   suivant qui, bien que dans le m\^eme esprit,  est sans lien direct avec les r\'esultats sur les connexions.

\begin{theoreme}
Sur les surfaces d'Inoue,  toutes les structures g\'eom\'etriques holomorphes sont localement
homog\`enes.
\end{theoreme}

\begin{demonstration}

Consid\'erons  une structure g\'eom\'etrique holomorphe $\tau$ sur une surface d'Inoue $S$. Soit  $\nabla_{0}$ l'unique  connexion   affine holomorphe sans torsion (et plate) sur  $S$. Nous
d\'emontrons que la structure g\'eom\'etrique $(\tau, \nabla_{0})$ qui consiste en la juxtaposition de $\tau$ et de $\nabla_{0}$ est localement homog\`ene. Autrement dit,
le pseudo-groupe des isom\'etries locales de $\nabla_{0}$ qui pr\'eservent $\tau$ agit transitivement sur $S$. L'avantage de ce proc\'ed\'e est de travailler avec  la structure g\'eom\'etrique $(\tau, \nabla_{0})$ qui est rigide au sens de Gromov~\cite{DG,Gro}. 

Par le th\'eor\`eme 3 de~\cite{D1},  les structures g\'eom\'etriques holomorphes rigides sur les vari\'et\'es complexes  compactes de dimension alg\'ebrique nulle  sont localement homog\`enes en dehors d'un sous-ensemble analytique compact d'int\'erieur vide (eventuellement vide). Il vient que $(\tau, \nabla_{0})$ est localement homog\`ene sur $S$ priv\'e
d'un sous-ensemble analytique compact d'int\'erieur vide  $E$.

Nous allons d\'emontrer que l'ensemble $E$ est vide.

Montrons d'abord  que, eventuellement \`a rev\^etement double non ramifi\'e  pr\`es,  $E$ est lisse (ceci est automatique si $S$ est constitu\'e d'un nombre fini de points, mais $S$ peut \'egalement admettre des composantes de dimension complexe un).

Supposons par l'absurde que $S$ n'est pas une sous-vari\'et\'e lisse de $E$.

Comme le pseudo-groupe des isom\'etries locales de $(\tau, \nabla_{0})$ pr\'eserve $E$, il vient que ce pseudo-groupe laisse invariant l'ensemble des points singuliers de $E$.

Consid\'erons $p \in E$ un point singulier de $E$. En particulier, $p$ n'est pas isol\'e dans $E$. Mais  $p$ est isol\'e parmi les points singuliers de $E$ et, par cons\'equent, $p$ est un point isol\'e dans son orbite sous l'action du pseudo-groupe des isom\'etries locales. Il d\'ecoule que 
 chaque champ de Killing local s'annule en $p$.

Notons $\mathcal G$ l'alg\`ebre des germes de champs de Killing au voisinage de $p$. Comme $\mathcal G$ agit transitivement sur un ouvert, elle est de dimension au moins $2$.

L'action  de $\mathcal G$ pr\'eserve  $\nabla$
  et  se lin\'earise donc en coordonn\'ees exponentielles au voisinage du point fixe  $p$. La lin\'earisation plonge $\mathcal G$ dans l'alg\`ebre de Lie de $GL(2, \CC)$. En particulier,
$\mathcal G$ est de dimension au plus $4$.

Si $\mathcal G$ est de dimension $2$, les sous-groupes de $GL(2, \CC)$ correspondants
sont conjugu\'es au  groupe des matrices diagonales, ou bien \`a l'un des  groupes suivants :

-  $\left(  \begin{array}{cc}
                                                                a   &   b \\
                                                                 0     &  a^{-1} \\
                                                                 \end{array} \right)$, avec $a \in \CC^*$ et $b \in \CC$;

- $\left(  \begin{array}{cc}
                                                                1   &   m \\
                                                                 0     &  n \\
                                                                 \end{array} \right)$, avec $m  \in \CC$ et $n \in \CC^*$;

 - $\left(  \begin{array}{cc}
                                                                m'   &   n' \\
                                                                 0     &  1 \\
                                                                 \end{array} \right)$, avec $m' \in \CC^*$ et $n' \in \CC$.

   Dans le premier cas,  le ferm\'e invariant $E$ s'identifie via l'application exponentielle  \`a la reunion des deux droites propres, tandis que dans
                                                                 le deuxi\`eme cas $E$ s'identifie \`a  la droite invariante $y=0$. Dans les deux situations, $E$ est lisse (quitte \`a consid\'erer un rev\^etement double
                                                                 non ramifi\'e de $S$).

     R\'eglons maintenant le cas o\`u $\mathcal G$ est de dimension $3$ ou $4$. Dans ce cas $\mathcal G$ engendre un sous-groupe de $GL(2, \CC)$ conjugu\'e au bien
     \`a $SL(2, \CC)$, ou bien \`a $GL(2, \CC)$ ou bien au groupe des matrices inversibles triangulaires sup\'erieures. Dans les deux premiers cas, il n'y a pas de ferm\'e
     invariant autre que le point $p$, qui doit par cons\'equent \^etre un point isol\'e (et donc lisse) de $E$ : absurde. Dans le dernier cas $E$ s'identifie comme avant \`a l'unique  droite invariante par les  matrices triangulaires sup\'erieures et est donc lisse.
     
     On vient de montrer que (\`a rev\^etement double pr\`es) $E$ est une sous-vari\'et\'e holomorphe de $S$. Si $E$ admet des composantes de dimension complexe un, alors
     ces composantes sont une reunion de courbes ferm\'ees (lisses). Or, les surfaces de Inoue ne contiennent aucune courbe ferm\'ee (lisse)~\cite{Inoue}.
     
     Il vient que $E$ est compos\'e d'un nombre fini de points. Supposons par l'absurde que $E$ est non vide et consid\'erons $p \in E$. On vient de montrer que
     l'alg\`ebre de Lie $\mathcal G$ est n\'ecessairement isomorphe \`a $SL(2, \CC)$ ou bien \`a $GL(2, \CC)$. En effet, dans tous les autres cas il existe au voisinage de $p$
     des droites  invariantes et donc $E$ ne se r\'eduit pas \`a un nombre fini de points.
     
     Traitons d'abord le cas o\`u $\mathcal G$ est l'alg\`ebre de Lie $sl(2, \CC)$ de $SL(2, \CC)$. L'action lin\'eaire de $\mathcal G$ sur $T_{p}S$ \'etant fid\`ele, elle s'identifie n\'ecessairement \`a l'action canonique de $sl(2, \CC)$
     sur $\CC^2$. Cette action a deux orbites : le point origine $p$ et $\CC^2 \setminus \{p  \}$. Le stabilisateur $H$ d'un \'el\'ement de $\CC^2 \setminus \{p  \}$ sous l'action correspondante  de $SL(2, \CC)$   est un sous-groupe \`a un param\`etre conjugu\'e dans
     $SL(2, \CC)$  \`a un sous-groupe unipotent de la forme $\left(  \begin{array}{cc}
                                                                1  &   b \\
                                                                 0     &  1 \\
                                                                 \end{array} \right)$, avec $b \in \CC$. 
                                                                 
            Remarquons que $G/H$ poss\`ede un champ de vecteurs holomorphe $G$-invariant qui s'exprime sur $\CC^2 \setminus \{p  \}$  sous la forme
               $\displaystyle z_{1}     \frac{\partial}{\partial z_{1} }  +  z_{2}  \frac{\partial}{\partial z_{2}}$.

           Comme l'ouvert $S\setminus E$ est localement model\'e sur $(G,G/H)$, il h\'erite d'un champ de vecteurs   holomorphe $X$. Le principe de prolongement
           de Hartog implique que le champ de vecteurs $X$ se prolonge \`a $S$. Par construction, $X$ est $\mathcal G$-invariant sur $S$ priv\'e de $E$ et, par analyticit\'e,
           $X$ doit \^etre  invariant partout. Ceci implique que $X$ s'annule au point $p$ car l'action de l'isotropie $SL(2, \CC)$ en $p$ ne pr\'eserve aucun vecteur non nul de $T_{p}S \simeq \CC^2$. La contradiction recherch\'ee vient du fait que les surfaces d'Inoue ne supportent  aucun champ de vecteurs non trivial qui s'annule en au moins un point~\cite{Inoue}.
           
           On peut \'egalement conclure en d\'emontrant directement que les surfaces d'Inoue n'admettent aucun feuilletage singulier $\mathcal F$ avec toutes les singularit\'ees locales de
           la forme $\displaystyle z_{1}     \frac{\partial}{\partial z_{1} }  +  z_{2}  \frac{\partial}{\partial z_{2}}$. En effet,  le th\'eor\`eme de Hartog permet de prolonger le fibr\'e tangent $T \mathcal F$ en un fibr\'e holomorphe en droites  sur $S$ (qui n'est pas un sous-fibr\'e de $TS$ \`a cause des singularit\'es)~\cite{GM} et, comme les classes de Chern d'une surface d'Inoue sont nulles~\cite{Barth-Peters, Wall}, les formules de Baum-Bott~\cite{GH} donnent $k=4k=c_{1}^2(T \mathcal F)$, o\`u $k$ est le cardinal de $S$ et $c_{1} (T \mathcal F)$ est la premi\`ere classe de Chern du  fibr\'e  $T \mathcal F$. Il  vient que $k=0$ et l'ensemble $S$ est vide.
           
           La preuve est  la m\^eme quand $\mathcal G$ est l'alg\`ebre de Lie de $GL(2, \CC)$.
\end{demonstration}

\section{ Fibr\'es elliptiques au-dessus d'une courbe de genre $g \geq 2$}

L'existence d'une structure  affine holomorphe  sur un fibr\'e principal holomorphe en courbes elliptiques au-dessus d'une surface de genre $g \geq 2$, de premier nombre
de Betti impair, est un r\'esultat d\^u \`a Maehara~\cite{Maehara}.
Dans~\cite{Klingler}  Klingler \'etudie  la g\'eom\'etrie globale des   structures   affines holomorphes  sur ces fibr\'es et montre, en particulier, 
que les applications d\'eveloppantes  associ\'ees  peuvent  \^etre non injectives :  le rev\^etement universel ne s'identifie pas toujours \`a un ouvert de $\CC^2$.

Nous \'etudions ici les connexions affines holomorphes non n\'ecessairement plates sur ces fibr\'es et d\'emontrons le point $(ii)$ du th\'eor\`eme~\ref{principal}.

Commen\c cons par suivre la construction g\'eom\'etrique de~\cite{Klingler} d'une connexion affine holomorphe sans torsion  plate sur $S$, qui nous servira de r\'ef\'erence par la suite.

\`A rev\^etement fini et quotient fini pr\`es, on peut supposer que le nombre  de Chern de $S$ vaut $g-1$~\cite{Klingler}. Alors le fibr\'e $S$ est  construit  de la mani\`ere suivante.
      
    Soit    $\Gamma$ un sous-groupe discret sans torsion de $PSL(2, \RR)$ tel que $\Sigma =\Gamma \backslash   H$, o\`u  $H$ d\'esigne le demi-plan de Poincar\'e. Choisissons une structure  projective holomorphe  sur  $\Sigma$, avec une application d\'eveloppante $\tau : H \to P^1(\CC)$ et un morphisme d'holonomie  $\rho : \Gamma \to PSL(2, \CC)$. L'image de $\Gamma$  dans  $PSL(2,\CC)$ se rel\`eve \`a  $SL(2, \CC)$ (ceci vient du fait que chaque vari\'et\'e (r\'eelle)  compacte orientable de dimension $3$ a une seconde 
    classe de Stiefel-Whitney triviale~\cite{MS}). 
    
    Consid\'erons donc $\Gamma$ comme \'etant un sous-groupe de $SL(2, \CC)$. D\'esignons par $W = \CC^2 \setminus \{0 \}$ le fibr\'e tautologique (de fibre $\CC^*$ ) au-dessus de  $P^1(\CC)$. L'action canonique de  $\Gamma$ sur $\CC^2$ induit une action de $\Gamma$ sur  $W$. On transporte la structure   affine holomorphe  de $W$ en  une structure affine   sur le fibr\'e  image r\'eciproque  $\tau^*(W) \simeq \CC^* \times H$.  La structure affine de $\tau^*(W)$ est invariante par l'action induite de $\Gamma$ (qui provient donc de l'action sur
    $H$ en tant que groupe de rev\^etement et de l'action sur $W$ don\'ee par $\rho$).
    
    Cette structure  affine holomorphe  sur  $\tau^*(W)$ est \'egalement invariante par les homoth\'eties dans les fibres.  Soit   un r\'eseau $\Delta \simeq \ZZ$ dans $\CC^*$ qui agit par multiplication dans les fibres de   $\tau^*(W)$
      et  consid\'erons  le quotient de  $\tau^*(W)$ by $\Delta \times \Gamma$. Ce quotient est un fibr\'e principal holomorphe au-dessus de $\Sigma$, de fibre la courbe elliptique $\Delta \backslash  \CC^*$, muni d'une structure affine holomorphe.
      
      \`A rev\^etement et quotient fini pr\`es, $S$ est biholomorphe au fibr\'e pr\'ec\'edent et h\'erite donc d'une connexion  affine holomorphe plate sans torsion $\nabla_{0}.$
      Le rev\^etement universel  $\CC \times H$ h\'erite \'egalement d'une  connexion qui provient de $\nabla_{0}$   via l'application
      
      $$\CC \times H  \to \CC^* \times H$$$$(z, \xi) \to (e^z, \xi).$$
      
      L'action d'un \'el\'ement  $\gamma=\left(  \begin{array}{cc}
                                                                a &   b \\
                                                                 c   &  d \\
                                                                 \end{array} \right)  \in SL(2, \RR)$ de $\Gamma$  sur  $\CC \times H$  est donn\'ee par la formule suivante~\cite{Klingler}:
                                                                 
      $$\gamma (z, \xi)=(z+log(c  \xi +d), \gamma \xi), \forall (z, \xi) \in \CC \times H,$$ o\`u $log$ d\'esigne  une  determination du  logarithme et l'action de $\gamma$ sur $H$ provient 
      de l'action canonique de $SL(2, \RR)$ sur $H$.
       
       Passons maintenant \`a la preuve du th\'eor\`eme~\ref{principal} :

\begin{demonstration}

La diff\'erence $\nabla - \nabla_{0}$  entre une connexion affine holomorphe arbitraire  $\nabla$ sur $S$ et la connexion de r\'ef\'erence $\nabla_{0}$ est une section holomorphe $\omega$ du fibr\'e vectoriel
$T^*S \otimes T^*S \otimes TS$. L'image r\'eciproque de ce tenseur sur le rev\^etement  universel $\CC  \times H$ de $S$ est une section du m\^eme type $\tilde{\omega}$ invariante
par le groupe de rev\^etement  engendr\'e par $\Delta \times \Gamma$ et par la transformation $t(z,\xi) = (z +2i \pi, \xi)$.

En les coordonn\'ees globales   $(z, \xi) \in \CC \times H$ du rev\^etement universel,  l'expression de $\tilde{\omega}$ sera  :
$$ f_{11}(z, \xi) dz \otimes dz \otimes \frac{\partial}{\partial z} + 
   f_{12}(z, \xi) dz \otimes d \xi  \otimes \frac{\partial}{\partial z}+
    f_{21}(z, \xi) dz \otimes dz \otimes \frac{\partial}{\partial \xi}+
     f_{22}(z, \xi) dz \otimes d \xi  \otimes \frac{\partial}{\partial \xi}+$$
    $$  g_{11}(z, \xi) d\xi  \otimes dz \otimes \frac{\partial}{\partial z}+
       g_{12}(z, \xi) d \xi  \otimes d \xi \otimes \frac{\partial}{\partial z}+
        g_{21}(z, \xi) d \xi  \otimes dz \otimes \frac{\partial}{\partial  \xi}+
         g_{22}(z, \xi) d \xi  \otimes d \xi  \otimes \frac{\partial}{\partial \xi},$$
         
         avec $f_{ij},g_{ij}$ des fonctions holomorphes quelque soient $1 \leq i,j \leq 2$.

        Par ailleurs,  un calcul direct montre que  la diff\'erence entre $\nabla_{0}$ et la connexion standard de $\CC \times H$ correspond aux fonctions 
            $f_{11}=f_{22}=g_{21}=1$, tous les autres coefficients \'etant nuls~\cite{Klingler}.
         
         L'invariance de $\tilde \omega$ par l'action de $\Delta$  implique que les fonctions 
         $f_{ij}(\cdot, \xi)$ and $g_{ij}(\cdot, \xi)$ sont invariantes par la transformation $z   \to z +log {\delta}$, o\`u $\delta$ est un \'el\'ement qui engendre $\Delta$. Comme ces fonctions doivent \^etre \'egalement invariantes par $t$, elles ont deux p\'eriodes  lin\'eairement ind\'ependantes sur $\RR$. 
        
         {\it  Le principe du maximum implique que 
        les fonctions  $f_{ij}$ and $g_{ij}$  ne d\'ependent que de la variable $\xi$}.
        
         Il vient que 
        $\frac{\partial}{\partial z}$  est un champ de Killing de $\tilde{\omega}$. Comme  le champ de vecteurs $\frac{\partial}{\partial z}$  pr\'eserve aussi  $\nabla_{0}$, il est  un champ de Killing  pour $\nabla$.

\`A partir de maintenant nous consid\'erons seulement des  connexions $\nabla$ sans torsion, ce qui conduit \`a : $f_{12}=g_{11}$ et  $f_{22}=g_{21}$.

   L'invariance de $\tilde{\omega}$ par l'action de $\Gamma$ est \'equivalente au syst\`eme  suivant :

                 $(1)   f_{11}(\xi)=f_{11}(\gamma \xi) -c f_{21}(\gamma \xi)(c \xi +d)$
                  
                $  (2)     f_{12}(\xi)=f_{12}(\gamma \xi) (c \xi +d)^{-2}-2c^2f_{21}(\gamma \xi)  -c f_{22} (\gamma \xi) (c \xi +d)^{-1} +2c f_{11}(\gamma \xi)(c \xi +d)^{-1}$
                  
                  $(3)     f_{21}(\xi)=f_{21}(\gamma \xi) (c \xi +d)^2$
                  
                 $ (4)     f_{22}(\xi)=2f_{21}(\gamma \xi)c (c \xi +d)  +f_{22}(\gamma \xi)$
                  
                 $ (5)     g_{12}(\xi)= g_{12}(\gamma \xi) (c \xi +d)^{-4} +c^2 f_{11}(\gamma \xi)(c \xi +d)^{-2}   +  c  f_{12}(\gamma \xi) (c \xi +d)^{-3}-c^3f_{21}(\gamma \xi) (c \xi +d)^{-1}  -c^2f_{22}(\gamma \xi)(c \xi +d)^{-2}  -cg_{22}(\gamma \xi)(c \xi +d)^{-3}$
                  
                  $(6) g_{22}(\xi)=g_{22}(\gamma \xi)(c \xi +d)^{-2} +cf_{22}(\gamma \xi)(c \xi +d)^{-1} +c^2f_{21}(\gamma \xi),$

   pour tout \'el\'ement  $\gamma =\left(  \begin{array}{cc}
                                                                a  &   b \\
                                                                 c   &  d\\
                                                                 \end{array} \right) \in \Gamma$.

           L'\'equation  $(3)$ montre que  $f_{21}$ a l'\'equivariance d'un champ de vecteurs holomorphe sur  $\Sigma$. Comme  $g \geq 2$, ceci implique que  $f_{21}=0$.
          Les \'equations $(1)$ et  $(4)$ impliquent que  $f_{11}$ et  $f_{22}$ sont des fonctions holomorphes sur $\Sigma$, donc constantes.
           
           Ceci conduit \`a :
           
           $(2')   f_{12}(\xi)=f_{12}(\gamma \xi) (c \xi +d)^{-2}  -c f_{22}  (c \xi +d)^{-1} +2c f_{11}(c \xi +d)^{-1}$
           
           $(5') g_{12}(\xi)= g_{12}(\gamma \xi) (c \xi +d)^{-4} +c^2 f_{11}(c \xi +d)^{-2}   +  c  f_{12}(\gamma \xi) (c \xi +d)^{-3}-c^2f_{22} (c \xi +d)^{-2}  -cg_{22}(\gamma \xi)(c \xi +d)^{-3}$
           
           $(6')  g_{22}(\xi)=g_{22}(\gamma \xi)(c \xi +d)^{-2} +cf_{22}(c \xi +d)^{-1}$, avec  $f_{11}, f_{22} \in \CC$.
           
           Nous montrons maintenant que $f_{11}=f_{22}=0$. Supposons par l'absurde que  $f_{22}  \neq 0$.
           
           La derni\`ere \'equation implique que  la $1$-forme diff\'erentielle $\displaystyle \alpha = (g_{22}(\xi) -\frac{f_{22}}{\xi- \bar{\xi}})d\xi$ est   $\Gamma$-invariante. Cette  $1$-forme descend donc sur $\Sigma$.
           
           Par ailleurs, $\displaystyle d\alpha=   \frac{f_{22}}{(\xi- \bar{\xi})^2}d\xi \wedge d\bar{\xi}.$ Il vient que $d \alpha$ est un multiple non nul de la forme volume hyperbolique de $\Sigma$ et donc 
            $\int_{\Sigma} d \alpha \neq 0$. Ceci est en contradiction avec la formule de  Stokes. 
           
           Il vient que  $f_{22}=0$. L'\' equation  $(2')$ implique de la m\^eme mani\`ere que $f_{11}=0$.
          
         L'interpr\'etation des \'equations  $(2')$ et  $(6')$ est \`a pr\'esent la suivante :  $f_{12}$ et  $g_{22}$ sont des sections du fibr\'e canonique $K_{\Sigma}$ de  $\Sigma$.
           
           Notons  $v=g_{22}-f_{12}$. L'\'equation  $(5')$ est \'equivalente au fait que   $w= -2 g_{12} + v'$  est une section holomorphe de $K_{\Sigma}^2$, autrement dit, $w(\xi)d\xi^2$ est  $\Gamma$-invariant.
           
           Il vient que l'espace des connexions affines holomorphes sans torsion sur $S$ est param\'etr\'e par $$(f_{12}, g_{22},-2g_{12}+g'_{22}-f'_{12})  \in 
            (H^0(\Sigma, K_{\Sigma}) )^2 \times H^0(\Sigma, K_{\Sigma}^{\otimes 2}) \simeq \CC^{5g-3}.$$

        Un calcul direct fournit le tenseur de courbure $R$ de la connexion : 
        
       $$R(\frac{\partial}{\partial z}, \frac{\partial}{\partial \xi}) \frac{\partial}{\partial z}=f_{12} \frac{\partial}{\partial z};~~R(\frac{\partial}{\partial z}, \frac{\partial}{\partial \xi}) \frac{\partial}{\partial \xi}=\lbrack f_{12} (g_{22}-f_{12})-f'_{12} \rbrack  \frac{\partial}{\partial z} -f_{12} \frac{\partial}{\partial \xi}.$$

          {\it La courbure de $\nabla$ est identiquement nulle si et seulement si $f_{12}=0$}. Dans ce cas on retrouve bien  la param\'etrisation de l'espace des  connexions
        plates donn\'ee  dans~\cite{Klingler}.

        {\it Il reste \`a montrer que l'alg\`ebre de Lie locale de $\nabla$ est engendr\'ee par $\frac{\partial}{\partial z}$, d\`es que $f_{12} \neq 0$.}
        
        Supposons donc que $f_{12} \neq 0$ et consid\'erons un champ de Killing local  $X=a(z, \xi) \frac{\partial}{\partial z}  + b(z, \xi) \frac{\partial}{\partial \xi}$, avec $a$ et $b$ des fonctions holomorphes locales sur $\CC \times H$.

  L'\' equation d'un champ de Killing pour $\nabla$  est~\cite{Kobayashi} :   $$\lbrack X, \nabla_{Y}Z\rbrack =\nabla_{\lbrack X,Y \rbrack}Z + \nabla_{Y}  \lbrack X, Z \rbrack,$$ quelque soient 
les champs de vecteurs tangents $Y,Z$. Il suffit de v\'erifier cette \'equation pour  $(Y,Z)$ \'etant un des couples suivants : $(\frac{\partial}{\partial z}, \frac{\partial}{\partial \xi}),
(\frac{\partial}{\partial z}, \frac{\partial}{\partial z})$ et  $(\frac{\partial}{\partial \xi}, \frac{\partial}{\partial \xi}).$

               L'\'equation aux d\'eriv\'ees partielles  satisfaite par $X=a(z, \xi) \frac{\partial}{\partial z}  + b(z, \xi) \frac{\partial}{\partial \xi}$ est  :

$ (i)~~a_{zz}  +  a_{z} + 2f_{12}b_{z}=0$

$(ii)~~  b_{z z}  + b_{z}=0$

$ \displaystyle (iii)~~a_{z \xi}    + g_{12}b_{z} +f_{12} b_{\xi}  + \frac{\partial f_{12}}{\partial  {\xi}}b=0$

$(iv)~~b_{z \xi}  +a_{z} + (g_{22}-f_{12})b_{z}=0$

$\displaystyle (v)~~a_{\xi \xi} -g_{12}a_{z}  +(2f_{12}-g_{12})a_{\xi}  +2g_{12} b_{\xi}  + \frac{\partial g_{12}}{\partial \xi}b=0$

$\displaystyle (vi)~~ b_{\xi \xi} + 2 a_{\xi}  +g_{22}b_{\xi} +  \frac{\partial g_{22}}{ \partial  \xi} b=0.$

La solution g\'en\'erale de l'\'equation $(ii)$ est   $b=\nu (\xi) e^{- z}  + C(\xi)$, o\`u $\nu, C$ sont des  fonctions holomorphes de la variable $\xi$.

On  remplace   $b_{z}$ dans l'\'equation $(i)$ et avec la m\'ethode de la variation de la constante on obtient    $a_{z}=  \lbrack z2f_{12}(\xi) \nu(\xi)+A(\xi) \rbrack e^{-z}$, o\`u $A $ est une fonction holomorphe de la variable $\xi$.

On introduit les d\'eriv\'ees partielles de $a$ et de  $b$ dans l'\'equation $(iv)$ :
$$  \lbrack - \nu'(\xi)  +A(\xi)   - ( g_{22} -2f_{12}) \nu (\xi)  \rbrack e^{-z}+ 2f_{12}(\xi) \nu(\xi)ze^{-z}=0.$$
 
 Comme les fonctions $e^{-z}$ et $ze^{-z}$  sont lin\'eairement ind\'ependantes sur $\CC$, ceci m\`ene \`a $f_{12}(\xi) \nu(\xi)=0$. Comme $f_{12} \neq 0$, il vient que $\nu(\xi)=0$.
 
 L'annulation du coefficient de $e^{-z}$ dans la relation pr\'ec\'edente implique alors $A(\xi)=0$.

On obtient donc  $b=C(\xi)$ et  $a=B(\xi)$, avec $B$ fonction holomorphe de la variable $\xi$.

Le syst\`eme initial  devient alors :

$(iii')~~f'_{12}b+f_{12}b'=0$

$(v')~a'' +(2f_{12}-g_{12})a'+   2g_{12}b'-g'_{12}b=0$

$(vi')~~b''  +2a'+g_{22}b'+g'_{22}b=0.$

La solution g\'en\'erale de  $(iii')$ est   $b=\frac{C_{0}}{f_{12}}$, avec $C_{0} \in \CC$.

Comme le genre  $g$ de  $\Sigma$ est sup\'erieur ou \'egal \`a  $2$,  toute section holomorphe de $K_{\Sigma}$  s'annule en au moins un point. Par cons\'equent,  $f_{12}$  s'annule en au moins un point $\xi_{0} \in H$. 

Pla\c cons nous dans un  voisinage ouvert d'un point  $(z,\xi_{0}) \in \CC \times H$. Les germes de champs de Killing (holomorphes)  en $(z, \xi_{0})$  satisfont tous $C_{0}=0$  et donc $b=0$. Il vient qu'au voisinage de ce point les germes de champs de Killing sont de la forme  
$B(\xi) \frac{\partial}{\partial z}$. D'apr\`es l'\'equation $(vi')$ on a que  la fonction $B$ est constante.

L'alg\`ebre des germes de champs de Killing au voisinage de $(z, \xi_{0})$ est donc de dimension un. Comme cette alg\`ebre est de dimension constante sur $S$~\cite{Amores, DG,  Gro}, sa dimension est (partout) \'egale \`a un.

L'alg\`ebre des champs de Killing est donc engendr\'ee
par le champ fondamental $\frac{\partial}{\partial z}$ de la fibration principale, d\`es que $f_{12} \neq 0$.
\end{demonstration}

Pour en d\'eduire  le corollaire~\ref{cconnexions} il suffit de prouver la proposition suivante :

\begin{proposition} (i) Les connexions affines holomorphes sans torsion sur un fibr\'e elliptique $S$ au-dessus d'une surface de Riemann de genre $g \geq 2$, de premier nombre de Betti impair,  sont projectivement plates. 

                                    (ii) Les connexions affines holomorphes sans torsion  invariantes par translations sur $\CC^2$  sont projectivement plates.
\end{proposition}

\begin{demonstration}

(i)  On reprend les notations de la preuve pr\'ec\'edente. La connexion projective associ\'ee \`a la connexion affine sans torsion $\nabla$ est
donn\'ee par l'\'equation diff\'erentielle de second ordre $\xi''=K^0(z, \xi) + K^1(z, \xi) \xi' + K^2(z, \xi)(\xi')^2+ K^3(z, \xi)(\xi')^3$, avec
$K^0=- f_{21}=0$, $K^1=f_{11}-2f_{22}$, $K^2=-(g_{22}-2f_{12})$ et $K^3=g_{12}$~\cite{Cartan}.

Il est classiquement connu d'apr\`es les travaux~\cite{Liouville},~\cite{Tresse},~\cite{Cartan} que cette connexion projective est plate si et seulement si
les invariants de Cartan sont identiquement nuls. Or, ces invariants qui valent $$L_{1}=2K^1_{z\xi}-K^2_{zz}-3K^0_{\xi \xi}-6K^0K^3_{z} -3K^3K^0_{z} + 3K^0K^2_{\xi}+3K^2K^0_{\xi}+
K^1K^2_{z}-2K^1K^1_{\xi},$$ $$L_{2}=2K^2_{z \xi}  -K^1_{\xi \xi}  -3K^3_{z z} + 6 K^3  K^0_{\xi}  + 3 K^{0}K^3_{\xi}-3K^3K^1_{z}  -3K^1K^3_{z}   -K^2K^1_{\xi} +2K^2K^2_{z}$$ sont
nuls car $K^0=K^1=0$  et $K^2, K^3$ sont des fonctions de $\xi$.

(ii) Dans ce cas les fonctions $K^0,K^1,K^2$ et $K^3$ sont constantes et on conclut comme au point pr\'ec\'edent.
\end{demonstration}

\section{Connexions projectives}   \label{Connexions projectives}

Rappelons que, d'apr\`es~\cite{KO1,KO2}, les surfaces complexes compactes munies d'une connexion projective holomorphe normale sont biholomorphes \`a : 
\`a $P^2(\CC)$, \`a une surface qui admet une connexion affine holomorphe ou bien \`a un quotient non ramifi\'e de l'espace  hyperbolique complexe $H^2_{\CC}$.

Commen\c cons  par  la proposition suivante qui nous sera utile par la suite :

\begin{proposition} \label{avecconnexion} Soit $S$ une surface complexe qui poss\`ede une connexion affine holomorphe. Alors toute connexion projective holomorphe normale sur $S$ est 
projectivement \'equivalente \`a une connexion affine holomorphe sans torsion sur $S$.
\end{proposition}

Le corollaire~\ref{cconnexions}  implique alors :

\begin{corollaire} \label{avecconnexionp} Si $S$ poss\`ede une connexion affine holomorphe, alors  toute connexion projective holomorphe normale sur $S$ est plate.
\end{corollaire}

\begin{demonstration}

Cette proposition est une cons\'equence directe des r\'esultats de~\cite{Gunning}. En effet, il est montr\'e dans~\cite{Gunning} que $S$ admet une connexion affine holomorphe si et
seulement si $dlog \Delta_{ij}$ est un \'el\'ement trivial dans $H^1(S, \Omega^1)$, o\`u $\Delta_{ij}$ est le $1$-cocycle du fibr\'e canonique \`a $S$ et $\Omega^1$ est
le faisceau des $1$-formes diff\'erentielles holomorphes (voir~\cite{Gunning}, page 96).

Par ailleurs, cette m\^eme condition suffit  pour  associer \`a toute  connexion projective holomorphe normale sur $S$, une connexion affine holomorphe sans torsion 
qui lui est projectivement \'equivalente (voir dans~\cite{KO1}, pages 78-79,  la formule explicite (3.6)). 
\end{demonstration}

D\'ecrivons maintenat  la connexion projective normale plate standard sur les quotients non ramifi\'e de $H^2_{\CC}$. Rappelons que l'espace hyperbolique complexe
$H^2_{\CC}=SU(2,1)/ S(U(2) \times U(1))$ s'identifie \`a l'ouvert $\{ z= \lbrack z_{0} : z_{1} \rbrack \in P^2(\CC),  |z_{0} |  +|z_{1} |  <  |z_{2} |  \}$ de $P^2(\CC)$ sur lequel
le groupe des biholomorphismes agit projectivement. Tout quotient compact de $H^2_{\CC}$ par un sous-groupe discret de $SU(2,1)$ h\'erite donc d'une  structure 
projective holomorphe  standard qui, de plus, est unique (\`a structure complexe fix\'ee)~\cite{Mok-Yeung} (voir \'egalement~\cite{KO3, Klingler2}). 

En rassemblant les arguments de~\cite{HM,McKay2,Mok-Yeung} (voir \'egalement~\cite{KO3,Klingler2}),  nous obtenons ici la proposition suivante, qui ensemble avec le corollaire~\ref{avecconnexionp},  ach\`eve la preuve du th\'eor\`eme~\ref{projectif}.

\begin{proposition} Soit $S$ une surface complexe compacte qui admet une connexion projective holomorphe normale, mais  ne poss\`ede aucune  connexion affine holomorphe. Alors 
$S$ est projectivement \'equivalente  \`a $P^2(\CC)$ muni de sa connexion projective plate standard ou bien \`a un quotient compact non ramifi\'e de $H^2_{\CC}$ muni
de sa connexion projective plate standard. 
\end{proposition}

\begin{demonstration}  D'apr\`es le r\'esultat de~\cite{KO1,KO2}, la surface $S$ est biholomorphe
\`a $P^2(\CC)$ ou bien \`a un quotient non ramifi\'e $H^2_{\CC}$.

-Si $S$ est biholomorphe \`a $P^2(\CC)$, alors le th\'eor\`eme principal de~\cite{HM} montre que $S$ admet une {\it unique} connexion projective holomorphe  normale qui est la connexion plate standard de $P^2(\CC)$. Un r\'esultat plus fort de McKay~\cite{McKay2} (qui suppose seulement l'existence d'une courbe rationnelle) implique que la seule connexion   de Cartan  projective holomorphe (et, en particulier,  l'unique connexion  projective holomorphe normale)   sur $P^2(\CC)$ est la connexion {\it plate} canonique.

-Si $S$ est un quotient de $H^2_{\CC}$, il est montr\'e dans~\cite{Mok-Yeung} (voir   \'egalement~\cite{KO3, Klingler2}) que $S$ admet une {\it unique} connexion    projective holomorphe.
La premi\`ere \'etape de la preuve consiste \`a interpr\'eter la diff\'erence entre deux connexions projectives, comme une section holomorphe du fibr\'e $Hom(L,S')$, o\`u $L$
est le fibr\'e en droites tautologique sur le projectivis\'e du fibr\'e tangent $\pi : PTS \to S$ et $S'$ est le sous-fibr\'e du fibr\'e  tangent holomorphe  \`a $PTS$ donn\'e par le noyau de la diff\'erentielle de $\pi$.
L'id\'ee est de consid\'erer la diff\'erence entre des g\'eod\'esiques projectives issues du m\^eme vecteur tangent \`a $S$ (voir~\cite{Mok-Yeung}, page 258, pour les d\'etails de la preuve).

La deuxi\`eme \'etape d\'emontre que le fibr\'e pr\'ec\'edent n'admet aucune section holomorphe non triviale~\cite{Mok-Yeung} (proposition 2.1) et aboutit \`a l'unicit\'e de la connexion projective
holomorphe.
\end{demonstration}

Mots-cl\'es : connexions affines holomorphes  -connexions projectives holomorphes- surfaces complexes- champs de Killing locaux.

Classification math. :  53B21, 53C56, 53A55.

\newpage

{\footnotesize

\vspace{2cm}

Sorin Dumitrescu

\rule{4cm}{.05mm}

D\'epartement de Math\'ematiques d'Orsay 

\'Equipe de Topologie et Dynamique

Bat. 425

\rule{4cm}{.05mm}

U.M.R.   8628  C.N.R.S.

\rule{4cm}{.05mm}

Univ. Paris-Sud (11)

91405 Orsay Cedex

France

\rule{4cm}{.05mm}

Sorin.Dumitrescu@math.u-psud.fr

 \end{document}